\documentclass[11pt]{article}
\usepackage{amssymb}
\usepackage{amsthm}
\usepackage{graphicx}
\usepackage{graphics}
\usepackage{amsthm}

\input pictex.tex
\newcommand{\esp}{\hspace{0.05cm}}
\newcommand{\s}{\sigma}

\theoremstyle{definition}

\newtheorem{thm}{Theorem}[section]

\newtheorem{prop}[thm]{Proposition}

\newtheorem{lem}[thm]{Lemma}

\newtheorem{rem}[thm]{Remark}

\newtheorem{defn}[thm]{Definition}

\newtheorem{cor}[thm]{Corollary}

\newtheorem{ex}[thm]{Example}

\newcommand{\vsp}{\vspace{0.1cm}}

\begin{document}

\date{}
\author{by Andr\'es Navas\footnote{Funded by the PBCT/Conicyt 
Research Network on Low Dimensional Dynamics.} \esp \esp \esp 
and \esp Crist\'obal Rivas, \esp \esp with an Appendix by Adam Clay}

\title{A new characterization of Conrad's property for group orderings, 
with applications}
\maketitle

\vspace{-0.87cm}

\begin{abstract} We provide a pure algebraic version of the dynamical 
characterization of Conrad's property given in \cite{order}. This approach allows 
dealing with general group actions on totally ordered spaces. As an application, 
we give a new and somehow constructive proof of a theorem first established by 
Linnell: an orderable group having infinitely many orderings has uncountably 
many. This proof is achieved by extending to uncountable orderable groups a 
result of the first author about orderings which may be approximated by their 
conjugates. This last result is illustrated by an example of an exotic 
ordering on the free group given by the third author in the Appendix. 
\end{abstract}

\vspace{0.34cm}

\noindent{\Large {\bf Introduction}}

\vspace{0.34cm}

In recent years, relevant progress has been made in the theory of (left) orderable 
groups. This has been achieved mainly by means of the use of a recently introduced 
mathematical object, namely the {\em space of group orderings} (see for instance 
\cite{dub,linnell,witte,sikora}). This space may be endowed with a natural topology  
(roughly, two orderings are close if they coincide over large finite sets), 
and the study of this topological structure should reveal some algebraic 
features of the underlying group. In \cite{order} it was realized that, 
for this study, the classical Conrad property for group orderings becomes 
relevant. Bringing ideas and techniques from the theory of codimension-one 
foliations, the `dynamical' insight of this property was revealed. 
Unfortunately, many proofs of \cite{order} 
are difficult to read for people with a pure algebraic view of 
orderable groups. More importantly, some of the results therein do not 
cover the case of {\em uncountable} groups. Indeed, the dynamical analysis   
of group orderings is done via the so-called `dynamical realization' 
of orderable groups as groups of homeomorphisms of the line, which 
is not available for general uncountable orderable groups.

Motivated by this, we develop here an algebraic counterpart of (part of) 
the analysis of \cite{order}. We begin by giving a new characterization of the 
Conrad property that is purely algebraic, 
although it has a dynamical flavor ({\em c.f.}, 
Theorem \ref{Conrad=noCrossing}). This leads naturally to the notion of Conradian 
actions on totally ordered spaces. A relevant example concerns the action of an 
ordered group on the space of cosets with respect to a convex subgroup. In this 
setting, we define the notion of Conradian extension ({\em c.f.}, Example 
\ref{convex-action}), and we generalize Conrad's classical theorem on 
the `level' structure of groups admitting Conradian orderings 
({\em c.f.,} Theorem \ref{super-conrad}, Corollary \ref{mas}). 

A relevant concept introduced in \cite{order} is the {\em Conradian soul}, which 
corresponds to the maximal subgroup of an ordered group that is convex and restricted 
to which the ordering is Conradian. In \cite{order}, a more geometrical view of this 
notion was given in the case of countable groups. Here we provide an analogous 
algebraic description which applies to general (possibly uncountable) ordered 
groups ({\em c.f.}, Theorem \ref{S=CS}). 

The Conradian soul was introduced as a main tool for dealing with the problem of 
approximating a group ordering by its conjugates. For instance, it was shown in 
\cite{order} that if the Conradian soul of an ordering on a non-trivial countable 
group is trivial, then this ordering is an accumulation point of its set of  
conjugates. The extension of this result to uncountable orderable groups appears 
here as Theorem \ref{primero}. We point out that an independent proof using 
completely different ideas was given by Adam Clay in \cite{clay}. 

Based on the work of Linnell \cite{linnell}, it was shown in \cite{order} that if an ordering 
on a group is isolated in the corresponding space of orderings, then its Conradian 
soul is `almost trivial', in the sense that it has only finitely many orderings. 
It is then natural to deal with ordered groups $(\Gamma,\preceq)$ for which the 
Conradian soul $C_{\preceq}(\Gamma)$ is non-trivial but has only finitely 
many orderings. If $\preceq$ is not Conradian, then 
to each of the orderings on $C_{\preceq}(\Gamma)$ 
corresponds an ordering on $\Gamma$ (roughly, the new orderings on $\Gamma$ 
are obtained by changing the original one on $C_{\preceq}(\Gamma)$ 
but preserving the set of elements bigger than the identity outside). As it was 
proved in \cite{order}, at least one of these orderings on $\Gamma$ is an accumulation 
point of its set of conjugates provided that $\Gamma$ is countable. Here we extend 
this result to the case of uncountable groups ({\em c.f.,} Corollary \ref{finalito}). 

The property of being approximated by its conjugates does not hold for all of the finitely 
many orderings on $\Gamma$ obtained by the preceding construction. A remarkable example 
illustrating this fact, namely the Dubrovina-Dubrovin ordering $\preceq_{DD}$ on braid 
groups $B_n$ \cite{dub}, was extensively studied from this point of view 
in \cite{order}. In the Appendix, Adam Clay provides a different kind of 
example, namely an `exotic' ordering $\preceq_C$ on the free group $F_2$. As it is 
the case of $\preceq_{DD}$, the Conradian soul of $\preceq_C$ is isomorphic to $\mathbb{Z}$, 
and $\preceq_C$ is not an accumulation point of the set of its conjugates. (This answers 
by the negative a question suggested in \cite[Remark 4.11]{order}.) The main difference 
between $\preceq_{DD}$ and $\preceq_C$ lies on the fact that $\preceq_{DD}$ is an isolated 
point of the (uncountable) space of orderings of $B_n$, while $\preceq_C$ is non-isolated 
in the (also uncountable) space of orderings of $F_2$. (Actually, the space of orderings 
of $F_2$ is homeomorphic to the Cantor set \cite{MC,order}.)

As a final application of our methods, we give a new proof of a theorem first 
established by Linnell \cite{linnell}: if a group has infinitely many orderings, 
then it has uncountably many. Linnell's proof uses an argument from General Topology 
for reducing the general case to that of Conradian orderings for which prior arguments 
by Zenkov \cite{zenkov} apply. To deal with the non Conradian case, we use our 
machinery on Conradian souls. Note that this was already done in \cite{order} 
for countable groups: Theorem \ref{linnell-general} here corresponds to the 
extension to the case of uncountable groups.


\section{Crossings and Conradian orderings}

\subsection{An equivalent for Conrad's property}
\label{nuevo-equivalente}

\hspace{0.45cm} Let $\preceq$ be an {\bf{\em ordering}} on a group $\Gamma$, that 
is, a total order relation which is invariant by left multiplication. Recall that 
$\preceq$ is said to be {\bf {\em Conradian}} if for all $f \succ 1$ and all 
$g \succ 1$ (for short, for all {\bf {\em positive elements}} $f,g$) there 
exists $n \!\in\! \mathbb{N}$ such that $fg^n \succ g$. (See however 
Remark \ref{n=2}.) A subgroup $\Gamma_0$ of $\Gamma$ is 
{\bf {\em $\preceq$-Conradian}} if the restriction
of $\preceq$ to it is a Conradian ordering. 

\vsp 

A {\bf {\em crossing}} for the ordered group $(\Gamma,\preceq)$ 
is a 5-uple $(f,g,u,v,w)$ of elements in $\Gamma$ such that:

\vsp

\noindent -- \esp $u \prec w \prec v$,

\vsp

\noindent -- \esp $g^n u \prec v$ \esp and \esp $f^n v \succ u$
\esp for every $n \in \mathbb{N}$,

\vsp

\noindent -- \esp there exist $M,N$ in $\mathbb{N}$ so that \esp
$f^N v \prec w \prec g^M u$.

\vspace{0.15cm}

\begin{rem} It follows from the third condition 
that neither $f$ nor $g$ can be equal to the identity.
\end{rem}

\begin{rem} If $(f,g,u,v,w)$ is a crossing, then the inequalities
$f^n v \succ u$ and $g^n u\prec v$ actually hold for {\em every} 
integer $n$. Indeed, we necessarily 
have $f v \prec v$, since in the other case we would have
$v \succ w \succ f^N v \succ f^{N-1} v \succ \ldots \succ f v \succ v,$
which is absurd. Therefore, for $n > 0$,
$$f^{-n} v \succ f^{n-1} v \succ \ldots f^{-1} v \succ v \succ u.$$
The inequality $g^{-n} u \prec v$ for $n > 0$ may be checked similarly.
\label{remark1}
\end{rem}

\begin{rem} The reason of the use of different type of letter for the elements $f,g$ 
and $u,v,w$ will become clear in \S \ref{conrad-general}. Somehow, $u,v,w$ should be 
thought of as `reference points' instead of genuine group elements (see Figure 1).
\end{rem}

\vspace{0.2cm}


\beginpicture

\setcoordinatesystem units <1cm,1cm>


\putrule from 1.5 -2.5 to 6.5 -2.5
\putrule from 1.5 2.5 to 6.5 2.5
\putrule from 1.5 -2.5 to 1.5 2.5
\putrule from 6.5 -2.5 to 6.5 2.5

\plot
1.5 0
1.625 0.05296
1.75 0.109
1.875 0.16056
2 0.216
2.125 0.26056
2.25 0.3009
2.375 0.35296
2.5 0.4 /

\plot
2.5 0.4
2.625 0.45296
2.75 0.509
2.875 0.56056
3 0.616
3.125 0.66056
3.25 0.7009
3.375 0.75296
3.5 0.8 /

\plot
3.5 0.8
3.625 0.85296
3.75 0.909
3.875 0.96056
4 1.016
4.125 1.06056
4.25 1.1009
4.375 1.15296
4.5 1.2 /

\plot
4.5 1.2
4.625 1.25296
4.75 1.309
4.875 1.36056
5 1.416
5.125 1.46056
5.25 1.5009
5.375 1.55296
5.5 1.6 
5.625 1.65296 /


\plot
6.5 -0.82519
6.375 -0.85296
6.125 -0.909
5.875 -0.96056
5.625 -1.016
5.375 -1.06056
5.125 -1.1009
4.875 -1.15296
4.625 -1.2 /

\plot
4.625 -1.2
4.375 -1.25296
4.125 -1.309
3.875 -1.36056
3.625 -1.416
3.375 -1.46056
3.125 -1.5009
2.875 -1.55296
2.625 -1.6
2.375 -1.65296 /


\setdots

\plot 1.5 -2.5
6.5 2.5 /

\putrule from 2.32 -1.68 to 5.68 -1.68
\putrule from 2.32 -1.68 to 2.32 1.68 
\putrule from 2.32 1.68 to 5.68 1.68  
\putrule from 5.68 -1.68 to 5.68 1.68

\put{Figure 1: A crossing} at 4 -3.5
\put{} at -4.2 0

\small


\put{$u$} at 1.5 -2.8
\put{$v$} at 6.5 -2.8
\put{$w$} at 4   -2.8
\put{$f^N v$} at 2.8 -2.8 
\put{$g^M u$} at 5.2 -2.8 
\put{$\bullet$} at 1.5 -2.5
\put{$\bullet$} at 6.5 -2.5 
\put{$\bullet$} at 4   -2.5 
\put{$\bullet$} at 2.8 -2.5 
\put{$\bullet$} at 5.2 -2.5 

\put{$f$} at 4.7 -0.9
\put{$g$} at 3.4 1

\endpicture


\vspace{0.5cm}

The next result is the natural analogue of both Propositions  
3.16 and 3.19 of \cite{order} in our setting.

\vspace{0.1cm}

\begin{thm} {\em The ordering $\preceq$ is Conradian
if and only if $(\Gamma,\preceq)$ admits no crossing.}
\label{Conrad=noCrossing}
\end{thm}

\begin{proof} Suppose that $\preceq$ is not Conradian, and
let $f,g$ be positive elements so that $fg^n \prec g$ for every $n
\in \mathbb{N}$. We claim that $(f,g,u,v,w)$ is a crossing for
$(\Gamma,\preceq)$ for the choice $u = 1$, and $v = f^{-1}g$, 
and $w = g^2$. Indeed:

\vsp

\noindent -- From $fg^2 \prec g$ one obtains
$g^2 \prec f^{-1} g$, and since $g \succ 1$, this gives \esp
$1 \prec g^2 \prec f^{-1}g,$  \esp that is, \esp $u \prec w \prec v$.

\vsp

\noindent -- From $f g^n \prec g$ one gets \esp $g^n \prec f^{-1} g$,
\esp that is, \esp $g^n u \prec v$ \esp (for every $n \in \mathbb{N}$);
on the other hand, since both $f,g$ are positive, we have $f^{n-1} g \succ 1$,
and thus \esp $f^n (f^{-1}g) \succ 1$, \esp that is, \esp $f^n v \succ u$ \esp
(for every $n \in \mathbb{N}$).

\vsp

\noindent -- The relation \esp $f (f^{-1}g) = g \prec g^2$ \esp may be read as
\esp $f^N v \prec w$ \esp for $N=1$; on the other hand, the relation \esp
$g^2 \prec g^3$ \esp is \esp \esp $w \prec g^M u$ \esp for $M = 3$.

\vsp

\vsp

Conversely, assume that $(f,g,u,v,w)$ is a crossing for $(\Gamma,\preceq)$ so that 
\esp $f^N v \prec w \prec g^M u$ \esp (with $M,N$ in $\mathbb{N}$). We will prove 
that $\preceq$ is not Conradian by showing that, for $h = g^Mf^N$ and $\bar{h} = g^M$, 
both elements $w^{-1} h w$ and $w^{-1} \bar{h} w$ are positive, but
$$(w^{-1} h w) (w^{-1} \bar{h} w)^n \prec w^{-1} \bar{h} w
\quad \mbox{ for all } \esp n \in \mathbb{N}.$$
To show this, first note that \esp $gw \succ w$. \esp 
Indeed, if not then we would have
$$w \prec g^N u \prec g^N w \prec g^{N-1} w \prec \ldots \prec gw \prec w,$$
which is absurd. Clearly, the inequality \esp $gw \succ w$ \esp implies
$$g^M w \succ g^{M-1} w \succ \ldots \succ gw \succ w,$$
and hence
\begin{equation}
w^{-1} \bar{h} w = w^{-1} g^M w \succ 1.
\label{zero}
\end{equation}
Moreover,
$$h w = g^M f^N w \succ g^M f^N f^N v 
= g^M f^{2N} v \succ g^{M} u \succ w.$$
and hence
\begin{equation}
w^{-1} h w \succ 1.
\label{one}
\end{equation}
Now note that, for every $n \in \mathbb{N}$,
$$h \bar{h}^n w = h g^{Mn} w \prec h g^{Mn} g^M u
= h g^{Mn+M} u \prec h v = g^M f^N v \prec g^M w = \bar{h} w.$$
After multiplying by the left by $w^{-1}$, the last inequality becomes 
$$(w^{-1} h w) (w^{-1} \bar{h} w)^n = 
w^{-1} h \bar{h}^n w \prec w^{-1} \bar{h} w,$$
as we wanted to check. Together with (\ref{zero}) and (\ref{one}),
this shows that $\preceq$ is not Conradian. 
\end{proof}

\vsp

\begin{rem} A fact that will be not used in this work is that, for every 
Conradian group ordering $\preceq$, one actually has $fg^2 \!\succ\! g$ 
for all positive elements $f,g$ ({\em i.e.,} one can take `$n\!=\!2$' in 
the original definition). The proof given in \cite[Proposition 3.7]{order} 
uses the fact that, if $f,g$ are positive elements for which $fg^2 \prec g$, 
then letting $h = fg$ one has $f h^n \prec h$ for all $n \in \mathbb{N}$. 
This is illustrated by Figure 2. Notice that, as shown below, in this 
situation $(f,fg,1,fg,g)$ is a crossing for $M=N=2$... 
\label{n=2} 
\end{rem}

\vspace{0.45cm}


\beginpicture

\setcoordinatesystem units <1.45cm,1.45cm>


\putrule from 1.5 -2.5 to 6.5 -2.5
\putrule from 1.5 2.5 to 6.5 2.5
\putrule from 1.5 -2.5 to 1.5 2.5
\putrule from 6.5 -2.5 to 6.5 2.5

\plot
1.5 1
1.625 1.05296
1.75 1.109
1.875 1.16056
2 1.216
2.125 1.26056
2.25 1.3009
2.375 1.35296
2.5 1.4 /

\plot
2.5 1.4
2.625 1.45296
2.75 1.509
2.875 1.56056
3 1.616
3.125 1.66056
3.25 1.7009
3.375 1.75296
3.5 1.8 /

\plot
3.5 1.8
3.625 1.85296
3.75 1.909
3.875 1.96056
4 2.016
4.125 2.06056 
4.25 2.1009
4.375 2.15296
4.5 2.2 /

\setquadratic
\plot
4.5 2.2   
4.75 2.329 
5 2.5 /

\plot 
1.5 -1.1 
3.8 -0.18 
5 0.2 /

\setlinear



\setquadratic
\plot 
1.5 -1.76  
5 -1.12   
6.5 0.2 /

\setlinear






\setdots

\plot 1.5 -2.5
6.5 2.5 /

\putrule from 5 -2.5 to 5 2.5  
\putrule from 1.5 1 to 6.5 1  
\putrule from 1.5 -1.12 to 6.5 -1.12    

\putrule from 2.35 -0.18 to 3.8 -0.18    
\putrule from 2.35 -0.18 to 2.35 -1.68
\putrule from 2.35 -1.68 to 3.8 -1.68 
\putrule from 3.8  -0.18 to 3.8 -1.68 
\putrule from 2.88 -2.5 to 2.88 -1.12 

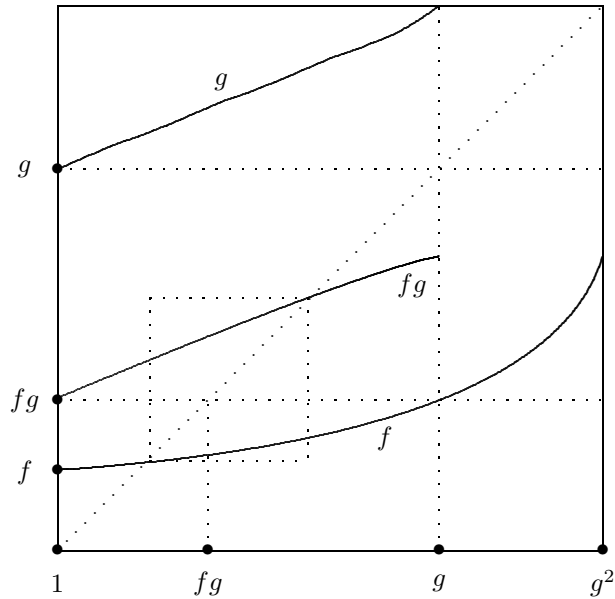
\put{Figure 2: The `$n\!=\!2$' condition} at 4 -3.5
\put{} at -1.47 0

\small


\put{$1$} at 1.5 -2.8 
\put{$g^2$} at 6.5 -2.8 
\put{$\bullet$} at 1.5 -2.5 
\put{$\bullet$} at 2.88 -2.5 

\put{$f$} at  1.2 -1.78    
\put{$fg$} at 1.2 -1.12 
\put{$fg$} at 2.88 -2.8 
\put{$g$} at  1.2 1  

\put{$fg$} at 4.75 -0.07   

\put{$g$} at 5 -2.8 
\put{$\bullet$} at 1.5 1
\put{$\bullet$} at 1.5 -1.12  
\put{$\bullet$} at 1.5 -1.76  
\put{$\bullet$} at 5 -2.5 
\put{$\bullet$} at 6.5 -2.5 

\put{$f$} at 4.5 -1.47 
\put{$g$} at 3 1.8

\endpicture


\vspace{0.25cm}

\begin{rem} The second condition in the definition of crossing may seem difficult to handle. 
A more `robust' property is that of {\bf {\em reinforced crossing}}, which is a 5-uple 
$(f,g,u,v,w)$ of elements in an ordered group $(\Gamma,\preceq)$ such that:

\vsp

\noindent -- \esp $u \prec w \prec v$,

\vsp

\noindent -- \esp $f u \succ u$ \esp and \esp $g(v) \prec v$,

\vsp

\noindent -- \esp there exist $M,N$ in $\mathbb{N}$ so that \esp
$f^N v \prec w \prec g^M u$.

\vspace{0.2cm}


\beginpicture

\setcoordinatesystem units <1cm,1cm>


\putrule from 1.5 -2.5 to 6.5 -2.5
\putrule from 1.5 2.5 to 6.5 2.5
\putrule from 1.5 -2.5 to 1.5 2.5
\putrule from 6.5 -2.5 to 6.5 2.5

\plot
1.5 0.2
6.5 1.96 /


\plot
6.5 0.2
1.5 -2.02 /


\setdots

\plot 1.5 -2.5
6.5 2.5 /

\putrule from 2.32 -1.68 to 5.68 -1.68
\putrule from 2.32 -1.68 to 2.32 1.68 
\putrule from 2.32 1.68 to 5.68 1.68  
\putrule from 5.68 -1.68 to 5.68 1.68

\put{Figure 3: A reinforced crossing} at 4 -3.5
\put{} at -4.2 0

\small


\put{$u$} at 1.5 -2.8
\put{$v$} at 6.5 -2.8
\put{$w$} at 4   -2.8
\put{$f^N v$} at 2.8 -2.8 
\put{$g^M u$} at 5.2 -2.8 
\put{$\bullet$} at 1.5 -2.5
\put{$\bullet$} at 6.5 -2.5 
\put{$\bullet$} at 4   -2.5 
\put{$\bullet$} at 2.8 -2.5 
\put{$\bullet$} at 5.2 -2.5 

\put{$f$} at 4.7 -0.9
\put{$g$} at 3.4 1.1

\endpicture


\vspace{0.45cm}

One easily checks that a reinforced crossing is a crossing. Conversely, if $(f,g,u,v,w)$ 
is a crossing, then $(f^{N} g^{M},g^{M} f^{N},f^N w, g^M w, w)$ is a reinforced 
crossing (here, $M,N$ in $\mathbb{N}$ are such that $f^N v \prec w \prec g^M u$). 
Indeed, from the properties of crossing one gets $f^{N} g^{M} (g^M w) \prec f^{N} v 
\prec w$ and $g^{M} f^{N} (f^N w) \succ g^{M} u \succ w$. Moreover, 
$f^{N} g^{M} (f^N w) \succ f^{N} g^{M} u \succ f^N w$ and 
$g^{M} f^{N} (g^M w) \prec g^{M} f^{N} v \prec g^M w$. 
\label{reinforced}
\end{rem}

\begin{rem} The dynamical characterization of Conrad's property may serve as 
inspiration for introducing other relevant properties for group orderings. 
(Compare \cite[Question 3.22]{order}.) 
For instance, one can say that a $6$-uple $(f,g,u,v,w_1,w_2)$ of 
elements in an ordered group $(\Gamma,\preceq)$ is a 
{\bf {\em (reinforced) double crossing}} if (see Figure 4):

\vsp

\noindent -- $u \prec w_1 \prec w_2 \prec v$, 

\vsp

\noindent -- $fu \succ u$ and $f v \succ v$,

\vsp

\noindent -- $gu \succ w_1$, \esp $gv \prec w_2$, \esp and $f w_2 \prec w_1$.

\vsp 
 
\noindent Finding a simpler algebraic counterpart of the property of not 
having a double crossing for an ordering seems to be an interesting problem. 
\end{rem}

\vspace{0.7cm}


\beginpicture

\setcoordinatesystem units <1cm,1cm>

\putrule from 0 0 to 6 0 
\putrule from 0 0 to 0 6 
\putrule from 6 0 to 6 6 
\putrule from 0 6 to 6 6

\put{$f$} at 4.6 2.5 
\put{$g$} at 1.5 3 

\put{$u$} at 0 -0.6
\put{$v$} at 6 -0.6 
\put{$w_1$} at 2.3 -0.6 
\put{$w_2$} at 3.7 -0.6

\put{$\bullet$} at 0 0 
\put{$\bullet$} at 6 0 
\put{$\bullet$} at 2.3 0 
\put{$\bullet$} at 3.7 0 

\put{Figure 4: A (reinforced) double crossing} at 3 -1.5 
\put{} at -5 -1 

\setquadratic

\plot 
0 2.5 
3 3 
6 3.4 /

\plot
0 0.5 
4 2.2 
6 6.8 /

\setlinear

\setdots

\plot 
3.7 0 
3.7 1.95 /

\plot 
1.95 1.95  
3.7 1.95 /

\plot 
1.95 1.95 
1.95 0 /
\plot 
0 0 
6 6 /

\plot 
6 6 
6 6.8 /

\plot 
2.5 0
2.5 2.5 /

\plot 
2.5 2.5 
0 2.5 /

\plot 
3.4 3.4 
6 3.4 / 

\plot 
3.4 3.4 
3.4 0 /

\endpicture


\subsection{An extension to group actions on ordered spaces}
\label{conrad-general}

\hspace{0.45cm} Let $\Gamma$ be a group acting by order-preserving 
bijections on a totally ordered space $(\Omega,\leq)$. A {\bf {\em crossing}} 
for the action of $\Gamma$ on $\Omega$ is a 5-uple $(f,g,u,v,w)$, where 
$f,g$ belong to $\Gamma$ and $u,v,w$ are in $\Omega$, such that:

\vsp

\noindent -- \esp $u \prec w \prec v$,

\vsp

\noindent -- \esp $g^n u \prec v$ \esp and \esp $f^n v \succ u$
\esp for every $n \in \mathbb{N}$,

\vsp

\noindent -- \esp there exist $M,N$ in $\mathbb{N}$ so that \esp
$f^N v \prec w \prec g^M u$.

\vsp

\begin{ex} The real line carries a natural total order, and thus our definition 
applies to groups acting on it by orientation preserving homeomorphisms. The 
notion of crossing for this case is exactly the same as that of {\em elements 
in transversal position} in \cite[Definition 3.24]{order}.
\end{ex}

\begin{ex} If $\Gamma$ is endowed with an ordering $\preceq$, one may take 
$(\Omega,\leq) = (\Gamma,\preceq)$ as a totally ordered set. The action of 
$\Gamma$ by left translations on it preserves the order: a crossing for 
this action corresponds to a crossing for $(\Gamma,\preceq)$, in the 
terminology of \S \ref{nuevo-equivalente}. Note that 
this example generalizes the preceding one for countable groups, 
since every countable ordered group may be canonically (up to semiconjugacy) realized 
as a group of orientation preserving homeomorphisms of the real line \cite[\S 2.1]{order}. 
\end{ex}

\vsp

For another relevant example recall that, given ordered group 
$(\Gamma,\preceq)$, a subset $S$ is {\bf {\em $\preceq$-convex}} if 
for every $f_1 \!\prec\! f_2 $ in $S$, every $f \!\in\! \Gamma$ satisfying 
\esp $f_1 \prec f  \prec f_2$ \esp belongs to $S$. When $S$ is a subgroup, 
this is equivalent to that for all positive $\bar{f} \!\in\! S$, every
$f \!\in\! \Gamma$ such that $1 \prec f \prec \bar{f}$ belongs to $S$. 

\vsp

\begin{ex} Let $(\Gamma,\preceq)$ be an ordered group, and let $\Gamma_0$ be a 
$\preceq$-convex subgroup. The space of left cosets $\Omega = \Gamma / \Gamma_0$ carries 
a natural total order $\leq$, namely $f \Gamma_0 < g \Gamma_0$ if $f h_1 \prec g h_2$ 
for some $h_1,h_2$ in $\Gamma_0$ (the reader will easily check that this definition is 
independent of the choice of $h_1$ and $h_2$ in $\Gamma_0$). The action of $\Gamma$ 
by left translations on $\Omega$ preserves this order. (Note that taking $\Gamma_0$ 
as being the trivial subgroup, this example reduces to the preceding one.) Whenever 
this action has no crossings, we will say that $\Gamma$ is a 
{\bf{\em $\preceq$-Conradian extension}} of $\Gamma_0$. 
\label{convex-action}
\end{ex}

\begin{rem} Let $(\Gamma,\preceq)$ be an ordered group, and let $\Gamma_0$ be 
a $\preceq$-convex subgroup. Given any ordering $\preceq_*$ on $\Gamma_0$, the 
{\bf{\em extension}} of $\preceq_*$ by $\preceq$ is the ordering $\preceq^*$ 
on $\Gamma$ for which $1 \prec^* f$ if and only if either $f \in \Gamma_0$ and 
$1 \prec_* f$, or $f \notin \Gamma_0$ and $1 \prec f$. The reader can easily 
check that $\Gamma_0$ is still a $\preceq^*$-convex subgroup of $\Gamma$. 
Moreover, $\Gamma$ is a $\preceq$-Conradian extension of $\Gamma_0$ if 
and only if it is a $\preceq^*$-Conradian extension of it.
\label{convex-extension}
\end{rem}

For a general order-preserving action of a group $\Gamma$ on a totally ordered space 
$(\Omega,\leq)$, the action of an element 
$f \!\in\! \Gamma$ is said to be {\bf {\em cofinal}} 
if for all $x < y$ in $\Omega$ there exists $n \!\in\! \mathbb{Z}$ such that 
$f^n(x) > y$. Note that if the action of $f$ is not cofinal, then there exist 
$x<y$ in $\Omega$ such that $f^n (x) < y$ for every integer $n$. 

\vspace{0.1cm}

\begin{prop} {\em Let $\Gamma$ be a group acting by order-preserving bijections on a totally 
ordered space $(\Omega,\leq)$. If the action of \esp $\Gamma$ on $\Omega$ has no crossings, 
then the set of elements whose action is not cofinal forms a normal subgroup of \esp $\Gamma$.}
\label{normality}
\end{prop}

\begin{proof} Let us denote the set of elements whose action is not cofinal 
by $\Gamma_0$. This set is normal. Indeed, given $g \in \Gamma_0$, let $x < y$ in 
$\Omega$ be such that $g^n(x) < y$ for all $n$. For each $h \in \Gamma$ we have 
$g^n h^{-1} (h(x)) < y$, and hence $(hgh^{-1})^n (h(x)) < h(y)$ (for all 
$n \!\in\! \mathbb{Z}$). Since $h(x) < h(y)$, this shows that $hgh^{-1}$ 
belongs to $\Gamma_0$.

It follows immediately from the definition that $\Gamma_0$ is stable under 
inversion, that is, $g^{-1}$ belongs to $\Gamma_0$ for all $g \!\in\! \Gamma_0$. 
The fact that $\Gamma_0$ is stable by multiplication is more subtle. For the 
proof, given $x \in \Omega$ and $g \in \Gamma_0$, we will denote by $I_g (x)$ 
the {\em convex closure} of the set $\{g^n(x)\!\!: n \in \mathbb{Z}\}$, that is, the 
set formed by the $y \in \Omega$ for which there exists $m,n$ in $\mathbb{Z}$ so 
that $g^m (x) \leq y \leq g^n (x)$. Note that $I_g(x) = I_g(x')$ for all 
$x' \in I_g (x)$; moreover, $I_{g^{-1}} (x) = I_g (x)$ for all $g \!\in\! \Gamma_0$ 
and all $x \!\in\! \Omega$; finally, if $g (x) = x$, then $I_g (x) = \{ x \}$. We 
claim that if $I_g (x)$ and $I_f (y)$ are non-disjoint for some $x,y$ in $\Omega$ 
and $f,g$ in $\Gamma_0$, then one of them contains the other. Indeed, 
assume that there exist non-disjoint sets $I_f (y)$ and $I_g (x)$, none of 
which contains the other. Without loss of generality, we may assume that $I_g (x)$ 
contains points to the left of $I_f (y)$ (if this is not the case, just interchange 
the roles of $f$ and $g$). Changing $f$ and/or $g$ by their inverses if necessary, 
we may assume that $g(x) >x$ and $f(y) < y$, and hence $g(x') > x'$ for all 
$x' \in I_g (x)$, and $f(y') < y'$ for all $y' \in I_y(f)$. Take 
$u \in I_g(x) \setminus I_f(y)$, $w \in I_g(x) \cap I_f(y)$, and 
$v \in I_f(y) \setminus I_g(x)$. Then one easily checks that 
$(f,g,u,v,w)$ is a crossing, which is a contradiction. 

Let now $g,h$ be elements in $\Gamma_0$, and let $x_1 < y_1$ and $x_2 < y_2$ be points in 
$\Omega$ such that $g^n (x_1) < y_1$ and $h^n (x_2) < y_2$ for all $n \!\in\! \mathbb{Z}$. 
Put $x = \min\{x_1,x_2\}$ and $y = \max\{y_1,y_2\}$. Then $g^n (x) < y$ and $h^n (x) < y$ 
for all $n \in \mathbb{Z}$; in particular, $y$ does not belong to neither $I_g (x)$ nor 
$I_h (x)$. Since $x$ belongs to both sets, we have either $I_g(x) \subset I_h(x)$ or 
$I_h(x) \subset I_g(x)$. Both cases being analogous, let us consider only the first 
one. Then for all $x' \in I_g(x)$ we have $I_h (x') \subset I_g (x') = I_g(x)$. 
In particular, $h^{\pm 1}(x')$ belongs to $I_g (x)$ for all $x' \in I_g (x)$. Since 
the same holds for $g^{\pm 1}(x')$, this easily implies that $(gh)^n(x) \in I_g(x)$ for 
all $n \in \mathbb{Z}$. As a consequence, $(gh)^n(x) < y$ for all $n \in \mathbb{Z}$, 
thus showing that $gh$ belongs to $\Gamma_0$. 
\end{proof}

\vspace{0.15cm}

Recall that for an ordered group $(\Gamma,\preceq)$, a {\bf{\em convex jump}} 
is a pair $(G,H)$ of distinct $\preceq$-convex subgroups such that $H$ is 
contained in $G$, and there is no $\preceq$-convex subgroup between them. 
The previously developed ideas lead naturally to the following result, 
which may be viewed as an extension of Conrad's theorem on the structure 
of convex subgroups for Conradian orderings \cite[Theorem 4.1]{conrad}. 
However, our proof follows ideas which are rather different from those 
of Conrad, and is much inspired from \cite[Exercise 2.2.46]{yo}.

\vspace{0.25cm}

\begin{thm} {\em Let $(\Gamma,\preceq)$ be an ordered group, and let $(G,H)$ 
be a convex jump in $\Gamma$. Suppose that $G$ is a Conradian extension of 
$H$. Then $H$ is normal in $G$, and the ordering induced by $\preceq$ on the 
quotient $G/H$ is Archimedean (and hence order isomorphic to a subgroup of 
$(\mathbb{R},+)$, due to H\"older's theorem} \cite{botto,kopytov,yo}{\em)}.
\label{super-conrad}
\end{thm}

\begin{proof} Let us consider the action of $G$ on the space of 
cosets $G /H$. Each element in $H$ fixes the coset $H$, and hence its action 
is not cofinal. By Proposition \ref{normality}, if we show that the action 
of each element in $G \setminus H$ is cofinal, then this will give the 
normality of $H$ in $G$. 

Now given $f \in G \setminus H$, let $G_f$ the smallest convex subgroup of 
$G$ containing $H$ and $f$. We claim that $G_f$ coincides with the set \esp 
$S_f \!=\! \{g \in G \!: f^m \prec g \prec f^n \mbox{ for some } 
m,n \mbox{ in } \mathbb{Z}\}.$ \esp 
Indeed, $S_f$ is clearly a convex subset of $G$ containing $H$ and contained in $G_f$. 
Thus, for showing that $G_f = S_f$, we need to show that $S_f$ is a subgroup. For this, 
first note that, in the notation of the proof of Proposition \ref{normality}, the 
conditions $g \in S_f$ and $I_g (H) \subset I_{f} (H)$ are equivalent. Therefore, 
for each $g \in S_f$ we have $I_{g^{-1}} (H) = I_g (H) \subset I_f (H)$, and thus 
$g^{-1} \! \in S_f$. Moreover, if $\bar{g}$ is another element in $S_f$, then 
$\bar{g} g H \in \bar{g} (I_f (H)) = I_f (H)$, and thus $I_{\bar{g}g} (H) 
\subset I_f (H)$. This means that $\bar{g} g$ belongs to $S_f$, thus 
concluding the proof that $S_f$ and $G_f$ coincide.

Each $f \in G \setminus H$ leads to a convex subgroup $G_f = S_f$ strictly containing $H$. 
Since $(G,H)$ is a convex jump, we necessarily have $S_f = G$. Given $g_1 \prec g_2$ in $G$, 
choose $m_1,n_2$ in $\mathbb{Z}$ for which $f^{m_1} \prec g_1$ and $g_2 \prec f^{n_2}$. 
Then we have $f^{n_2-m_1} g_1 \succ f^{n_2-m_1} f^{m_1} = f^{n_2} \succ g_2$, and hence 
$f^{n_2 - m_1} (g_1 H) \geq g_2 H$. This easily implies that the action of 
$f$ is cofinal.

We have then show that $H$ is normal in $G$. The left invariant total order on the 
space of cosets $G / H$ is therefore a group ordering. Moreover, given $f,g$ in $G$, 
with $f \notin H$, the previous argument shows that there exists $n \!\in\! \mathbb{Z}$ 
such that $f^n \succ g$, and thus $f^n H \succeq gH$. This is nothing but the 
Archimedean property for the induced ordering on $G / H$. 
\end{proof}

\vsp 

\begin{cor} {\em Under the hypothesis of Theorem} \ref{super-conrad}, {\em up to 
multiplication by a positive real number, there exists a unique nontrivial group 
homomorphism $\tau \!: G \rightarrow \mathbb{R}$ such that $ker(\tau) \!=\! H$ 
and $\tau(g) \!>\! 0$ for every positive element $g \in G \setminus H$.}
\label{mas}
\end{cor}


\section{On the approximation of a group ordering by its conjugates}

\subsection{Describing the Conradian soul via crossings}

\hspace{0.45cm} The {\bf{\em Conradian soul}} \esp 
$C_{\preceq}(\Gamma)$ of an ordered group $(\Gamma,\preceq)$ 
corresponds to the (unique) subgroup which is $\preceq$-convex, $\preceq$-Conradian, 
and which is maximal among subgroups verifying these two properties simultaneously. 
This notion was introduced in \cite{order}, where a dynamical counterpart in the case 
of countable groups was given. To give an analogous characterization in the general 
case, we consider the set $S^+$ formed by the elements $h \!\succ\! 1$ 
such that $h \preceq w$ for every crossing $(f,g,u,v,w)$ satisfying  
$1 \preceq u$. Analogously, we let $S^{-}$ be the set formed by the 
elements $h \prec 1$ such that $w \preceq h$ for every crossing
$(f,g,u,v,w)$ satisfying $v \preceq 1$. Finally, we let  
$$S = \{ 1 \} \cup S^{+} \cup S^{-}.$$   
{\em A priori}, it is not clear that the set $S$ has a nice structure 
(for instance, it is not at all evident that it is actually a 
subgroup). However, this is largely shown by the theorem below.

\vspace{0.1cm}

\begin{thm} {\em The Conradian soul of $(\Gamma,\preceq)$ coincides
with the set $S$ above.}
\label{S=CS}
\end{thm}

\vspace{0.1cm}

Before passing to the proof, we give four general lemmas on crossings for 
group orderings (note that the first three lemmas still apply to crossings 
for actions on totally ordered spaces). The first one allows us replacing the 
`comparison'element $w$ by its `images' under positive iterates of either 
$f$ or $g$.

\vspace{0.1cm}

\begin{lem} {\em If $(f,g,u,v,w)$ is a crossing,
then $(f,g,u,v,g^n w)$ and $(f,g,u,v,f^n w)$ are
also crossings for every $n \!\in\! \mathbb{N}$.}
\label{lema1}
\end{lem}

\begin{proof} We will only consider the first 5-uple (the 
case of the second one is analogous). Recalling that $g w \succ w$,
for every $n \! \in \! \mathbb{N}$ we have $u \prec w \prec g^n
w$; moreover, $v \succ g^{M+n} u = g^n g^M u \succ g^n w$. Hence,
\esp $u \prec g^n w \prec v$. \esp On the other hand, $f^N v \prec w 
\prec g^n w$, while from $g^M u \succ w$ we get $g^{M+n} u \succ g^n w$. 
\end{proof}

\vspace{0.3cm}

Our second lemma allows replacing the `limiting' elements $u$ and $v$
by more appropriate ones.

\vspace{0.15cm}

\begin{lem} {\em Let $(f,g,u,v,w)$ be a crossing. If $f u \succ u$ (resp. 
$f u \prec u$) then $(f,g,f^n u,v, w)$ (resp. $(f,g,f^{-n} u, v, w)$) is also 
a crossing for every $n > 0$. Analogously, if $g v \prec v$ (resp. $g v \succ v$), 
then $(f, g, u, g^n v, w)$ (resp. $(f,g,u,g^{-n}v,w)$) is also crossing for every 
$n > 0$.} 
\label{lema3}
\end{lem}

\begin{proof} Let us only consider the first 5-uple (the case of the 
second one is analogous).
Suppose that $fu \succ u$ (the case $fu \prec u$ may be treated similarly).
Then $f^n u \succ u$, which gives $g^M f^n u \succ g^M u \succ w$. To show
that $f^n u \prec w$, assume by contradiction that $f^n u \succeq w$. Then
$f^n u \succ f^N v$, which gives $u \succ f^{N-n} v$, which is absurd. 
\end{proof}

\vspace{0.3cm}

The third lemma relies on the dynamical insight of the crossing condition. 

\vspace{0.15cm}

\begin{lem}{\em If $(f,g,u,v,w)$ is a crossing, then $(hfh^{-1},
hgh^{-1},hu,hv,hw)$ is also a crossing for every $h \in \Gamma$.}
\label{lema2}
\end{lem}

\begin{proof} The three conditions to be checked are nothing but
the three conditions in the definition of crossing multiplied by $h$
by the left. 
\end{proof}

\vspace{0.3cm}

A direct application of the lemma above shows that, if $(f,g,u,v,w)$ is a crossing, then the
5-uples $(f,f^ngf^{-n}, f^n u, f^n v , f^n w)$ and $(g^nfg^{-n}, g ,g^n u, g^n v, g^n w)$
are also crossings for every $n \in \mathbb{N}$. This combined with Lemma \ref{lema3}
may be used to show the following.

\vspace{0.15cm}

\begin{lem} {\em If $(f,g,u,v,w)$ is a crossing and $1 \preceq h_1 \prec h_2$ are
elements in $\Gamma$ such that $h_1 \in S$ and $h_2 \notin S$, then there exists a
crossing $(\tilde{f},\tilde{g},\tilde{u},\tilde{v},\tilde{w})$ such that
$h_1 \prec \tilde{u} \prec \tilde{v} \prec h_2$.} \label{lemapro}
\end{lem}

\begin{proof} Since $1\prec h_2  \notin S$, there must be a
crossing $(f,g,u,v,w)$ such that $1\preceq u \prec w \prec h_2$. Let
$N \in \mathbb{N}$ be such that $f^N v\prec w $. Denote by
$(f,\bar{g},\bar{u},\bar{v},\bar{w})$ the crossing
$(f, f^N g f^{-N}, f^N u, f^N v, f^N w)$. Note that
$\bar{v} = f^N v \prec w \prec h_2$.
We claim that $h_1 \preceq \bar{w} = f^N w$. Indeed, if $f^N u \succ u$ then
$f^n u \succ 1$, and by the definition of $S$ we must have $h_1 \preceq \bar{w}$.
If $f^N u \prec u$, then we must have $f u \prec u$, so by Lemma
\ref{lema3} we know that $(f,\bar{g},u,\bar{v},\bar{w})$ is also a 
crossing, which allows still concluding that $h_1 \preceq \bar{w}$.

Now for the crossing $(f,\bar{g},\bar{u}, \bar{v}, \bar{w})$ there exists
$M \in \mathbb{N}$ such that $\bar{w} \prec \bar{g}^M \bar{u}$. Let us
consider the crossing $(\bar{g}^M f \bar{g}^{-M}, \bar{g}, \bar{g}^M \bar{u},
\bar{g}^M \bar{v}, \bar{g}^M \bar{w})$. If $\bar{g}^M \bar{v} \prec \bar{v}$
then $\bar{g}^M \bar{v} \prec h_2$, and we are done. If not, then
we must have $\bar{g} \bar{v} \succ \bar{v}$. By Lemma \ref{lema3},
$(\bar{g}^M f \bar{g}^{-M}, \bar{g}, \bar{g}^M \bar{u}, \bar{g}^M \bar{v}, \bar{w})$ 
is still a crossing, and since $\bar{v} \prec h_2$, this concludes the proof. 
\end{proof}


\vspace{0.38cm}

\noindent{\em Proof of Theorem} \ref{S=CS}. The proof is divided into several steps.

\vsp\vsp

\noindent {\underbar{Claim 0.}} The set $S$ is convex.

\vsp

This follows directly from the definition of $S$.

\vsp\vsp

\noindent {\underbar{Claim 1.}} If $h$ belongs
to $S$, then $h^{-1}$ also belongs to $S$.

\vsp 

Assume that $h \in S$ is positive and $h^{-1}$ does not belong to
$S$. Then there exists a crossing $(f,g,u,v,w)$ so that $h^{-1}
\prec w \prec v \preceq 1$.

We first note that, if $h^{-1} \preceq u$, then after conjugating by 
$h$ as in Lemma \ref{lema2}, we get a contradiction because $(hgh^{-1},
hfh^{-1}, hu, hv, hw)$ is a crossing with \esp $1 \preceq hu $ \esp and
\esp $hw \prec hv \preceq h$. To reduce the case $h^{-1} \succ u$ to this
one, we first
use Lemma \ref{lema2} and we consider the crossing $(g^Mfg^{-M}, g, g^M u,
g^M v, g^M w)$. Since \esp $h^{-1} \prec w \prec g^M u \prec g^M w \prec g^M v$,
\esp if $g^M v \prec v$ then we are done. If not, Lemma \ref{lema3} shows that
$(g^Mfg^{-M}, g, g^M u , g^M v, w)$ is also a crossing, which still
allows concluding.

For the case where $h \in S$ is negative 
({\em i.e.,} its inverse is positive) we proceed similarly but we
conjugate by $f^N$ instead of $g^M$. Alternatively, since $1 \in S$ 
and $1\prec h^{-1}$, if we suppose that $h^{-1}\notin S$ then Lemma
\ref{lemapro} provides us with a crossing $(f,g,u,v,w)$ such
that $1\prec u\prec w \prec v \prec h^{-1}$, which gives a
contradiction after conjugating by $h$.

\vsp\vsp

\noindent {\underbar{Claim 2.}} If $h$ and $\bar{h}$ 
belong to $S$, then $h\bar{h}$ also belongs to $S$.

\vsp 

First we show that for every positive elements in $S$, 
their product still belongs to $S$. (Note that, by Claim 1, 
the same will be true for products of negative elements in $S$.) 
Indeed, suppose that $h,\bar{h}$ are positive elements, with $h \in S$
but $ h \bar{h} \notin S$.
Then, by Lemma \ref{lemapro} we may produce a crossing $(f,g,u,v,w)$
such that $h \prec u \prec v \prec h \bar{h}$.
After conjugating by $h^{-1}$ we obtain the crossing
$(h^{-1}fh,h^{-1}gh,h^{-1}u,h^{-1}v,h^{-1}w)$
satisfying $1 \prec h^{-1}u \prec h^{-1} w \prec \bar{h}$, 
which shows that $\bar{h} \notin S$. 

Now, if $h \prec 1 \prec \bar{h}$ then 
$h \prec h\bar{h}$. Hence, if $h\bar{h}$ is negative then
the convexity of $S$ gives $h\bar{h} \in S$. If $h\bar{h}$ 
is positive, then $\bar{h}^{-1}h^{-1}$
is negative, and since $\bar{h}^{-1} \prec \bar{h}^{-1} h^{-1}$, 
the convexity gives again that $\bar{h}^{-1}h^{-1}$, and hence 
$h\bar{h}$, belongs to $S$. The remaining case 
$\bar{h} \prec 1 \prec h$ may 
be treated similarly.

\vsp \vsp

\noindent {\underbar{Claim 3.}} The subgroup $S$ is Conradian.

\vsp 

In order to apply Theorem \ref{Conrad=noCrossing}, we need to show that there are no
crossings in $S$. Suppose by contradiction that $(f,g,u,v,w)$ is a crossing such that 
$f,g,u,v,w$ all belong to $S$. If $1 \preceq w$ then, by Lemma \ref{lema2}, we have 
that $(g^n f g^{-n}, g, g^n u, g^n v, g^n w)$ is a crossing. Taking $n = M$ so 
that $g^M u \succ w$, this gives a contradiction to the definition of $S$ 
because $1 \preceq w \prec g^M u \prec g^M w \prec g^M v \in S$.
The case $w \preceq 1$ may be treated in an analogous way 
by conjugating by powers of $f$ instead of $g$.

\vsp\vsp

\noindent {\underbar{Claim 4.}} The subgroup $S$ is maximal 
among $\preceq$-convex, $\preceq$-Conradian subgroups.

\vsp 

Indeed, if $C$ is a subgroup strictly containing $S$, then there 
is a positive $h$ in $C \setminus S$. By Lemma \ref{lemapro},
there exists a crossing $(f,g,u,v,w)$ such that $1\prec u \prec w
\prec v \prec h$. If $C$ is convex, then $u,v,w$ belong to $C$.
To conclude that $C$ is not Conradian, it suffices to show
that $f$ and $g$ belong to $C$.

Since $1\prec u$, we have either $1 \prec g \prec g u \prec v$ or
$1 \prec g^{-1} \prec g^{-1} u \prec v$. In both cases, the convexity of
$C$ implies that $g$ belongs to $C$. On the other hand, if $f$ is positive  
then from $f^N \prec f^N v \prec w$ we get $f \in C$, whereas in the case
of a negative $f$ the inequality $1\prec u$ gives $1\prec f^{-1} \prec f^{-1}
u \prec v$, which still shows that $f \in C$. $\hfill \square$


\subsection{Approximation of group orderings: the role of the Conradian soul}

\hspace{0.45cm} For a (left) orderable group $\Gamma$, we denote by $\mathcal{LO}(\Gamma)$  
the set of all orderings on $\Gamma$. This space carries the topology having as a 
subbasis the family of sets $U_f \!=\! \{\preceq  :\esp f \succ 1\}$, where 
$f \neq 1$. Endowed with this topology, $\mathcal{LO}(\Gamma)$ is called 
the {\bf {\em space of (left) orderings}} of the group $\Gamma$.

\begin{rem} As shown in \cite{yo}, a simple application of Tychonov's theorem shows 
that $\mathcal{LO}(\Gamma)$ is always compact. Moreover, the $`n\!=\!2$' property from 
Remark \ref{n=2} implies that the subset of Conradian orderings is closed therein 
(and hence compact). A more dynamical argument for showing this consists in noticing 
that the condition that $(f,g,u,v,w)$ is a reinforced crossing for prescribed $M,N$ 
is clearly open in $\mathcal{LO}(\Gamma)$ ({\em c.f.}, Remark \ref{reinforced}).
\end{rem}

The {\bf {\em positive cone}} of an ordering $\preceq$ in $\mathcal{LO}(\Gamma)$ is the 
set $P$ of its positive elements. Because of the left invariance, $P$ completely 
determines $\preceq$. The {\bf {\em conjugate}} of $\preceq$ by $h \in \Gamma$ 
is the ordering $\preceq_h$ having positive cone $hPh^{-1}$. In other words, 
$g \succ_h 1$ holds if and only if $h g h^{-1} \succ 1$. We will say that $\preceq$ may 
be approximated by its conjugates if it is an accumulation point of its set of conjugates.

\vspace{0.2cm}

\begin{thm} {\em If the Conradian soul of an ordered group $(\Gamma,\preceq)$ is trivial 
and $\preceq$ is not Conradian, then $\preceq$ may be approximated by its conjugates.}
\label{primero}
\end{thm}

\begin{proof} Let $f_1 \prec f_2 \prec \ldots \prec f_k$ be finitely many 
positive elements in $\Gamma$. We need to show that there exists a conjugate of $\preceq$
which is different from $\preceq$ but for which all the $f_i$'s are still positive. Since 
$1 \!\in\! C_{\preceq}(\Gamma)$ and $f_1 \notin C_{\preceq}(\Gamma)$, Theorem \ref{S=CS} 
and Lemma \ref{lemapro} imply that there is a crossing $(f,g,u,v,w)$ such that
$1\prec u \prec v \prec f_1$. Let $M,N$ in $\mathbb{N}$ be such
that $ f^N v \prec w \prec g^M u$. We claim that $1\prec_{v^{-1}} f_i$
and $1\prec_{w^{-1}} f_i$ for $1\leq i \leq k$, but $g^M f^N \prec_{v^{-1}} 1$
and $g^M f^N \succ_{w^{-1}} 1$. Indeed, 
since $1 \prec v \prec f_i$, we have $v\prec f_i\prec f_i v$, thus $1
\prec v^{-1} f_i v$. By definition, this means that $f_i\succ_{v^{-1}} 1$.
The inequality $f_i \succ_{w^{-1}} 1$ is proved similarly. Now note that
$g^M f^N v \prec g^M w \prec v$, and so $g^M f^N \prec_{v^{-1}} 1$. Finally, from
$g^Mf^N w \succ g^M u \succ w $ we get $g^Mf^N \succ_{w^{-1}} 1$. 

Now the preceding relations imply that the $f_i$'s are still
positive for both $\preceq_{v^{-1}}$ and $\preceq_{w^{-1}}$, 
but at least one of these orderings is different from 
$\preceq$. This concludes the proof. 
\end{proof}

\vspace{0.3cm}

Based on the work of Linnell \cite{linnell}, it is shown in \cite[Proposition 4.1]{order} 
that no Conradian ordering is an isolated point of the space of orderings of a group 
having infinitely many orderings. Together with Theorem \ref{primero}, this 
shows the next proposition by means of the convex extension procedure 
({\em c.f.}, Remark \ref{convex-extension}). 

\vsp

\begin{prop} {\em Let $\Gamma$ be an orderable group. If $\preceq$ 
is an isolated point of $\mathcal{LO} (\Gamma)$, then its Conradian 
soul is non-trivial and has only finitely many orderings.}
\label{linelito}
\end{prop}

\vsp 

As a consequence of a nice theorem of Tararin, the number of orderings on an orderable 
group having only finitely many orderings is a power of 2; moreover, all of these 
orderings are necessarily Conradian \cite{kopytov,yo}. By the preceding theorem, if 
$\preceq$ is an isolated point of an space of orderings $\mathcal{LO} (\Gamma)$, 
then its Conradian soul admits $2^n$ different orderings for some $n \geq 1$, all 
of them Conradian. Let $\{\preceq_1,\preceq_2,\ldots,\preceq_{2^n}\}$ be these 
orderings, where $\preceq_1$ is the restriction of $\preceq$ to its Conradian soul. 
Since $C_{\preceq}(\Gamma)$ is $\preceq$-convex, each $\preceq_j$ induces an ordering 
$\preceq^j$ on $\Gamma$, namely the convex extension of $\preceq_j$ by $\preceq$. 
(Note that $\preceq^1$ coincides with $\preceq$.) All the orderings $\preceq^j$ 
share the same Conradian soul \cite[Lemma 3.37]{order}. Assume throughout that 
$\preceq$ is not Conradian.

\vsp\vsp

\begin{thm} {\em  With the notation above, at least one of the orderings
$\preceq^j$ is an accumulation point of the set of conjugates of $\preceq$.}
\label{final}
\end{thm}

\vsp

\begin{cor} {\em At least one of the orderings 
$\preceq^j$ is approximated by its conjugates.}
\label{finalito}
\end{cor}

\begin{proof} Asumming Theorem \ref{final}, we have 
$\preceq^k \in\! acc (orb(\preceq^1))$ for some $k \in \{1,\ldots,2^n\}$. 
Theorem \ref{final} applied to this $\preceq^k$ instead of $\preceq$ shows 
the existence of $k' \in \{1,\ldots,2^n\}$ so that 
$\preceq^{k'} \in acc (orb(\preceq^k))$, and hence
$\preceq^{k'} \in acc(orb(\preceq^1))$. If $k'$ equals either $1$ or $k$ then we are
done; if not, we continue arguing in this way... In at most $2^n$ steps we will find 
an index $j$ such that $\preceq^j \in acc(orb(\preceq^j))$. 
\end{proof}

\vspace{0.2cm}

Theorem \ref{final} will follow from the next proposition.

\vsp

\begin{prop} {\em Given an arbitrary finite family 
$\mathcal{G}$ of $\preceq$-positive elements 
in $\Gamma$, there exists $h \in \Gamma$ and 
$1 \prec \bar{h} \notin C_{\preceq} (\Gamma)$ such 
that $1 \prec h^{-1} f h \notin C_{\preceq} (\Gamma)$ 
for all $f \in \mathcal{G} \setminus C_{\preceq} (\Gamma)$, 
but $1 \succ h^{-1}\bar{h}h \notin C_{\preceq} (\Gamma)$.} 
\label{a-probar}
\end{prop}

\vsp

\noindent{\em Proof of Theorem} \ref{final} {\em from Proposition} 
\ref{a-probar}. Let us consider the directed set formed by the 
finite sets $\mathcal{G}$ of $\preceq$-positive elements. For each 
such a $\mathcal{G}$, let $h_{\mathcal{G}}$ and $\bar{h}_{\mathcal{G}}$
be the elements in $\Gamma$ provided by Proposition \ref{a-probar}. After 
passing to subnets of $(h_{\mathcal{G}})$ and $(\bar{h}_{\mathcal{G}})$ 
if necessary, we may assume that the restrictions of  
$\preceq_{h_{\mathcal{G}}^{-1}}$ to $C_{\preceq}(\Gamma)$ 
all coincide with a single $\preceq_j$. Now the properties 
of $h_{\mathcal{G}}$ and $\bar{h}_{\mathcal{G}}$ imply:

\vsp

\noindent -- $f \succ^j 1$ \esp \esp and \esp \esp 
$f \esp (\succ^j)_{h_{\mathcal{G}}^{-1}} \esp 1$ \esp 
\esp for all \esp $f \in \mathcal{G} \setminus C_{\preceq} (\Gamma)$,

\vsp

\noindent -- $\bar{h}_{\mathcal{G}} \succ 1$, \esp \esp  but \esp \esp \esp 
$\bar{h}_{\mathcal{G}} \esp \esp (\prec^j)_{h_{\mathcal{G}}^{-1}} \prec 1$.

\vsp

\noindent This clearly shows the Theorem. $\hfill\square$

\vspace{0.5cm}

For the proof of Proposition \ref{a-probar} we will use three general lemmas.

\begin{lem} {\em For every $1\prec c \notin C_{\preceq} (\Gamma)$ 
there is a crossing $(f,g,u,v,w)$ such that $u,v,w$ do not belong 
to $C_{\preceq}(\Gamma)$ and $1\prec u\prec w \prec v \prec c$.}
\label{primer-lema}
\end{lem}

\begin{proof} By Theorem \ref{S=CS} and Lemma \ref{lemapro}, 
for every $1\preceq s\in C_{\preceq} (\Gamma)$ there exists a 
crossing $(f,g,u,v,w)$ such that $s\prec u\prec w\prec v \prec c$. 
Clearly, $v$ does not belong to $ C_{\preceq} (\Gamma)$. The element 
$w$ is also ouside $ C_{\preceq} (\Gamma)$, since in the other case the 
element $a = w^2$ would satisfy $w \prec a \in C_{\preceq} (\Gamma)$, 
which is absurd. Taking $M > 0$ so that $g^M u \succ w$, this gives \esp 
$g^M u \notin C_{\preceq} (\Gamma)$, \esp $g^M w \notin C_{\preceq} (\Gamma)$, 
\esp and $g^M v \notin C_{\preceq} (\Gamma)$. \esp Consider 
the crossing $(g^Mfg^{-M}, g, g^M u,g^M v, g^M w)$. If $g^M v\prec v$, then 
we are done. If not, then $gv \succ v$, and Lemma \ref{lema3} ensures that 
$(g^Mfg^{-M}, g, g^M u,  v, g^Mw)$ is also a crossing, which still allows 
concluding. 
\end{proof}


\begin{lem} {\em Given $1 \prec c \notin C_{\preceq} (\Gamma)$ there exists 
$1 \prec a \notin C_{\preceq}(\Gamma)$ (with $a \prec c$) such that, for all 
$1 \preceq b \preceq a$ and all \esp $\bar{c} \succeq c$, one has 
$1 \prec b^{-1} \bar{c} b \notin C_{\preceq} (\Gamma)$.}
\label{segundo-lema}
\end{lem}

\begin{proof} Let us consider the crossing $(f,g,u,v,w)$ 
such that $1\prec u\prec w \prec v \prec c$ and such that $u,v,w$ 
do not belong to $ C_{\preceq} (\Gamma)$. We affirm that the Lemma holds 
for $a = u$ (actually, it holds for $a = w$, but the proof is slightly more
complicated). Indeed, if $1 \preceq b \preceq u$, then from $b
\preceq u \prec v \prec \bar{c}$ we get $1 \preceq b^{-1}u \prec
b^{-1}v \prec b^{-1} \bar{c}$, and thus the crossing
$(b^{-1}fb,b^{-1}gb,b^{-1}u,b^{-1}v,b^{-1}w)$ shows that $b^{-1}
\bar{c} \notin C_{\preceq} (\Gamma)$. Since $1 \preceq b$, we conclude that $1 \prec
b^{-1} \bar{c} \preceq b^{-1} \bar{c} b$, and the convexity of $ S$
implies that $b^{-1} \bar{c} b \notin C_{\preceq} (\Gamma)$. 
\end{proof}


\begin{lem} {\em For every $g\in \Gamma$ the set $g \esp C_{\preceq} (\Gamma)$ 
is convex. Moreover, for every crossing $(f,g,u,v,w)$ one has 
$u C_{\preceq} (\Gamma) < w C_{\preceq} (\Gamma) < v C_{\preceq} (\Gamma)$, 
in the sense that \esp $uh_1 \prec wh_2 \prec vh_3$ \esp for all $h_1,h_2,h_3$ 
in $ C_{\preceq} (\Gamma)$ (c.f., Example {\em \ref{convex-action}}).}
\label{lema4}\end{lem}

\begin{proof} The verification 
of the convexity of $g C_{\preceq} (\Gamma)$ is straightforward. Now
suppose that $uh_1 \succ wh_2$ for some $h_1,h_2$ in $ C_{\preceq} (\Gamma)$. 
Then since $u \prec w$, the convexity of both left classes 
$u C_{\preceq} (\Gamma)$ and $w C_{\preceq} (\Gamma)$ gives 
the equality between them. In particular, there exists 
$h \in C_{\preceq} (\Gamma)$ such that $uh = w$. Note that such 
an $h$ must be positive, so that $1 \prec h = u^{-1} w$. But since
$(u^{-1}fu, u^{-1}gu, 1,u^{-1}v, u^{-1}w)$ is a crossing, this 
contradicts the definition of $ C_{\preceq} (\Gamma)$. Showing that 
$w C_{\preceq} (\Gamma) \prec v C_{\preceq} (\Gamma)$ is similar. 
\end{proof}

\vspace{0.15cm}

\noindent{\em Proof of Proposition} \ref{a-probar}. Let us label the elements of 
$\mathcal{G} \!=\! \{f_1,\ldots,f_r\}$ so that $f_1 \prec \ldots \prec f_r$, 
and let $k$ be such that $f_{k-1} \in C_{\preceq} (\Gamma)$ but 
$f_{k} \notin C_{\preceq} (\Gamma)$. Recall that, 
by Lemma \ref{segundo-lema}, there exists $1 \prec a \notin C_{\preceq} (\Gamma)$ such
that, for every $1 \preceq b \preceq a$, one has $1 \prec b^{-1}
f_{k+j} b \notin C_{\preceq} (\Gamma)$ for all $j \geq 0$. 
We fix a crossing $(f,g,u,v,w)$ such that 
$1\prec u\prec v\prec a$ and $u\notin C_{\preceq} (\Gamma)$. Note that the conjugacy 
by $w^{-1}$ gives the crossing $(w^{-1}fw, w^{-1}gw,w^{-1}u,w^{-1}v,1)$.

\vsp\vsp

\noindent{\underbar{Case 1.}} One has $w^{-1}v \preceq a$. 

In this case, the proposition holds for \esp 
$h=w^{-1}v$ \esp and \esp $\bar{h}=w^{-1}g^{M+1}f^N w$. \esp To show this, 
first note than neither $w^{-1}gw$ nor $w^{-1}fw$ belong to $C_{\preceq} (\Gamma)$. 
Indeed, this follows from the convexity of $C_{\preceq} (\Gamma)$ and the 
inequalities \esp $w^{-1}g^{-M}w \prec w^{-1}u \notin C_{\preceq} (\Gamma)$ 
\esp and $w^{-1} f^{-N} w \succ w^{-1} v \notin C_{\preceq} (\Gamma)$. 
\esp We also have $1\prec w^{-1}g^{M}f^{N}w$, and hence $1\prec
w^{-1}gw \prec w^{-1} g^{M+1}f^{N}w$, which shows that 
$\bar{h} \notin C_{\preceq} (\Gamma)$. On the other hand, the 
inequality $w^{-1} g^{M+1}f^{N}w (w^{-1} v) \prec w^{-1}v$ reads as
$h^{-1}\bar{h}h\prec 1$. Finally, Lemma \ref{lema1} applied to the 
crossing $(w^{-1}fw, w^{-1}gw,w^{-1}u,w^{-1}v,1)$ shows that 
$(w^{-1} fw, w^{-1}gw, w^{-1} u, w^{-1} v, w^{-1} g^{M+n}f^{N}w)$ is 
a crossing for every $n > 0$. For $n \geq M$ we have $w^{-1}
g^{M+1}f^{N}w (w^{-1}v)\prec w^{-1} g^{M+n}f^{N}w $. Since 
$w^{-1} g^{M+n}f^{N}w \prec w^{-1} v$, Lemma \ref{lema4} easily 
implies that $w^{-1} g^{M+1}f^{N}w (w^{-1}v) C_{\preceq} (\Gamma) \prec 
w^{-1}vC_{\preceq} (\Gamma)$, that is, $h^{-1}\bar{h}h\notin C_{\preceq} (\Gamma)$.

\vsp\vsp

\noindent{\underbar{Case 2.}} One has $a\prec w^{-1}v$, 
but $w^{-1} g^m w \preceq a$ for all $m > 0$. 

We claim that, in this case, the proposition holds for $h = a$ and
$\bar{h}=w^{-1}g^{M+1}f^N w$. This may be checked in the very same
way as in Case 1 by noticing that, if $a\prec w^{-1}v$ but $w^{-1}g^m
w \succeq a$ for all $m > 0$, then $(w^{-1}fw, w^{-1}gw,w^{-1}u,a,1)$ 
is a crossing.

\vsp\vsp

\noindent{\underbar{Case 3.}} One has $a\prec w^{-1}v$ 
and $w^{-1} g^m w\succ a$ for some $m > 0$. (Note 
that the first condition follows from the second one.) 

We claim that, in this case, the proposition holds for $h=a$ and 
$\bar{h}=w \notin C_{\preceq} (\Gamma)$.
Indeed, we have $g^{m}w\succ ha$ (and $w\prec ha$), and
since $g^{m}w\prec v\prec a$, we have $wa\prec a$, which means that
$h^{-1}\bar{h}h \prec 1$. Finally, from Lemmas \ref{lema1} and \ref{lema4} we get \esp 
$wa C_{\preceq} (\Gamma)\preceq g^m w C_{\preceq} (\Gamma) 
\prec v C_{\preceq} (\Gamma) \preceq a C_{\preceq} (\Gamma)$. \esp This implies that \esp 
$a^{-1} w a C_{\preceq} (\Gamma) \prec C_{\preceq} (\Gamma)$, \esp which means that 
$h^{-1} \bar{h} h \notin C_{\preceq} (\Gamma)$. $\hfill\square$


\section{Finitely many or uncountably many group orderings}

\hspace{0.45cm} The goal of this final short section is to use 
the previously developed ideas to show the following result.

\vspace{0.2cm}

\begin{thm} {\em If the space of orderings of an 
orderable group is infinite, then it is uncountable.}
\label{linnell-general}
\end{thm}

\begin{proof} Let us fix an ordering $\preceq$ on an 
orderable group $\Gamma$. We need to analize two different cases. 

\vspace{0.2cm}

\noindent{\underbar{Case 1.}} The Conradian soul of $C_{\preceq}(\Gamma)$ 
is non-trivial and has infinitely many orderings.

\vspace{0.1cm}

This case was settled in \cite{order} (see Proposition 4.1 therein) using 
ideas going back to Zenkov \cite{zenkov} and Tararin \cite{kopytov}. 

\vspace{0.2cm}

\noindent{\underbar{Case 2.}} The Conradian soul of $C_{\preceq}(\Gamma)$ 
has only finitely many orderings.

\vspace{0.1cm}

If $\preceq$ is Conradian, then $\Gamma = C_{\preceq} (\Gamma)$ has finitely many 
orderings. If not, then Theorems \ref{primero} and \ref{final} imply that there exists 
an ordering $\preceq^*$ on $\Gamma$ which is an accumulation point of its conjugates. 
The closure in $\mathcal{LO}(\Gamma)$ of the set of conjugates of $\preceq^*$ is 
then a compact set without isolated points. By a well-known fact in General 
Topology, such a set must be uncountable. Therefore, $\Gamma$ admits 
uncountably many orderings. 
\end{proof}


\begin{small}

\vspace{0.1cm}


\noindent Andr\'es Navas\\

\noindent Dep. de Matem\'aticas, Fac. de Ciencia, Univ. de Santiago, 
Alameda 3363, Est. Central, Santiago, Chile\\

\noindent Email address: andres.navas@usach.cl\\

\vspace{0.37cm}

\noindent Crist\'obal Rivas\\

\noindent Dep. de Matem\'aticas, Fac. de Ciencias, Univ. de Chile, 
Las Palmeras 3425, \~Nu\~noa, Santiago, Chile\\

\noindent Email address: cristobalrivas@u.uchile.cl

\end{small}


\vspace{1cm}

\appendix

\section{\hspace{-0.8cm} 
ppendix: An exotic ordering of the free group on two elements,}

\vspace{-0.24cm}

${}$ \hspace{2.6cm} {\Large \bf by Adam Clay}

\begin{abstract}
We construct an ordering of $F_2$ which is not an accumulation point of its conjugates 
in $\mathcal{LO}(F_2)$ and whose Conradian soul is isomorphic to $\mathbb{Z}$.  
This ordering is realized as the restriction of the Dubrovina-Dubrovin 
ordering of $B_3$ to an appropriate free subgroup of $B_3$.
\end{abstract}

\vsp\vsp

We begin by defining the Dehornoy ordering of the braid groups (also known as the `standard' ordering), 
whose positive cone we shall denote $P_D$ \cite{PD94,DDRW02}. Recall that for each integer $n \ge 2$, 
the Artin braid group $B_n$ is the group generated by 
$\s_1 , \s_2 , \dots , \s_{n-1}$, subject to the relations 
$$\s_i\s_j = \s_j\s_i {\rm \: \esp if \esp \:} |i-j| >1,
\quad \s_i\s_j\s_i=\s_j\s_i\s_j {\rm \: \esp if \esp \: } |i-j| =1.$$  

\begin{defn} Let $w$ be a word in the generators $\s_i, \cdots , \s_{n-1}$ (so no $\s_j$ occurs 
for $j < i$). Then $w$ is said to be: 
$i$-positive if the generator $\s_i$ occurs in $w$ with only positive exponents, $i$-negative if 
$\s_i$ occurs with only negative exponents, and $i$-neutral if $\s_i$ does not occur in $w$.
\end{defn}

We then define the positive cone of the Dehornoy ordering as

\begin{defn} 
The positive cone $P_D \subset B_n$ of the Dehornoy ordering is the set 
\[P_D = \{ \beta \in B_n \!: \mbox{ $\beta$ is $i$-positive for some $i \leq n-1$}\}.\]
\end{defn}

An extremely important property of this ordering is that the conjugate $\beta \s_k \beta^{-1}$ is always 
$i$-positive for some $i$, for every generator $\s_k$ in $B_n$ and any braid $\beta \in B_n$. This 
property is referred to as the subword property \cite{DDRW02}.

There is also a second ordering of interest, discovered by the authors of \cite{dubr}, whose positive 
cone we shall denote by $P_{DD}$. Denote by $P_i \subset B_n$ the set of all $i$-positive braids.  
Note that the set of all $i$-negative braids is simply $P_i^{-1}$.

\begin{defn} The positive cone $P_{DD} \subset B_n$ is the set
\[P_{DD} = P_{1} \cup P_{2}^{-1} \cup \cdots \cup P_{n-1}^{(-1)^{n}}.\]
\end{defn}

From \cite{CP01}, we know that the subgroup of $B_3$ generated by the elements $\s_1^2, \s_2^2$ is isomorphic to $F_2$, the free group on two generators.  Thus we may consider $F_2$ to be the subgroup of $B_3$ generated by $\s_1^2$ and $\s_2^2$, and define a positive cone $P$ in $F_2$ by $P = P_{DD} \cap F_2$.  Note that any element of $F_2$ must always be represented by a braid word having even total exponent, and that the ordering $\preceq_C$ of $F_2$ asssociated to the positive cone $P$ is simply 
the restriction of the $P_{DD}$ ordering to the subgroup $\langle \s_1^2, \s_2^2 \rangle$.

\begin{prop}
{\em The ordering $\preceq_C$ is not an accumulation point of its conjugates in $\mathcal{LO}(F_2)$. 
Specifically, no conjugates of $\preceq_C$ distinct from $\preceq_C$ lie inside the open set 
$U_{\s_2^{-2}} \subset \mathcal{LO}(F_2)$.}
\end{prop}

\begin{proof} 
Let $\beta \in F_2 \subset B_3$ be given, and consider the positive cone $\beta P \beta^{-1}$. To prove 
the claim, we must show  that $\s_2^{-2} \in \beta P \beta^{-1}$ implies  $\beta P \beta^{-1} = P$.  

First, observe that conjugation of $P$ by any even power of $\s_2$ does not change $P$:  this follows 
from the fact that $\s_2^{-2}$ is the least positive element in the associated ordering $\preceq_C$ of 
$F_2$. Indeed, for any element $g \in P$, we have $\s_2^{-2} \preceq_C g$, so that $\s_2^2g \succeq_C 1$, 
and hence $\s_2^2 g \s_2^{-2} \succ_C 1$, that is, $\s_2^2 g \s_2^{-2} \in P$.

Now with $\s_2^{-2} \in \beta P \beta^{-1}$, in particular we must have 
$\beta^{-1} \s_2^{-2} \beta \in P$.  Since $P$ consists of those elements of $F_2$ that are either $1$-positive or $2$-negative, by the subword property, we know that  $\beta^{-1} \s_2^{-2} \beta$ is not $1$-positive, and so must be $2$-negative.  
Therefore  $\beta^{-1} \s_2^{-2} \beta = \s_2^{k}$ for some $k<0$, and in fact, by considering 
total exponents we see that  $\beta^{-1} \s_2^{-2} \beta = \s_2^{-2}$.

Recall that we are working in $F_2 \subset B_3$, so $\beta$ cannot commute with $\s_2^{-2}$ unless $\beta$ itself 
is an even power of $\s_2$ (the power must be even since $ \beta \in F_2 = \langle \s_1^2, \s_2^2 \rangle$).  
Therefore  $\s_2^{-2} \in \beta P \beta^{-1}$ implies $\beta = \s_2^{2k}$, so that $\beta P \beta^{-1} = P$.
\end{proof}

Next, we show that the only non-trivial convex subgroup in the ordering $\preceq_C$ 
of $F_2$ defined by $P$ is $\langle \s_2^{-2} \rangle$, the infinite cyclic group 
generated by the least positive element $\s_2^{-2}$.  In particular, this shows that 
the Conradian soul of the ordering $\preceq_C$ of $F_2$ is isomorphic to $\mathbb{Z}$.

\begin{thm} 
{\em Suppose that $S$ is a subgroup of $F_2$ that is convex in the ordering $\preceq_C$. 
If $S$ properly contains $\langle \s_2^{-2} \rangle$, then $S = F_2$.} 
\end{thm}

\begin{proof}
Let $S$ be a convex subgroup properly containing  $\langle \s_2^{-2} \rangle$.  As the containment is proper, 
$S$ must contain a $1$-positive braid $\beta$.  Suppose that $\beta$ is represented by the $1$-positive braid word $\s_2^k \s_1 w$
where $k \in \mathbb{Z}$, and $w$ is a $1$-positive, $1$-neutral or empty braid word.  Left multiplying by an appropriate power of $\s_2^2$, we may produce a new $1$-positive braid $\beta' = \s_2^{2l} \beta$ in $S$ that is represented by a $1$-positive braid word of the form $\s_2^{k'} \s_1 w$, where $k' = 2l + k >0$. 
Note that $\beta' \in S$, since both $\beta$ and $\s_2^2$ lie in $S$.  

Consider the braid represented by the word $\s_1^{-2}\s_2^{k'} \s_1 w$.  For any $m$, we have $\s_1^{-1} \s_2^m \s_1 = \s_2 \s_1^m \s_2^{-1}$, so that we compute
\[\s_1^{-2}\s_2^{k'} \s_1 w = \s_1^{-1}  \s_2 \s_1^{k'} \s_2^{-1} w =   \s_1^{-1}  \s_2 \s_1 \s_1^{k'-1} \s_2^{-1} w = \s_2 \s_1 \s_2^{-1} \s_1^{k'-1} \s_2^{-1} w,\]
and since $k'>0$ and $w$ is a $1$-positive, $1$-neutral or empty word, we see that $\s_1^{-2}\s_2^{k'} \s_1 w$ represents a $1$-positive braid.  Therefore, in the ordering $\preceq_C$ of $F_2$, we have
\[ 1 \prec_C \s_1^{-2}\s_2^{k'} \s_1 w  \Rightarrow  \s_1^2 \prec_C \s_2^{k'} \s_1 w = \beta' .\]
Since $1 \prec_C \s_1^2$ and $\beta' \in S$, we conclude that $\s_1^2 \in S$.  
But then $S$ contains both $\s_2^2$ and $\s_1^2$, the generators of $F_2$, so that $S = F_2$.
\end{proof}

\begin{rem}
From the work of \cite{MC2,order2}, it is known that the ordering $\preceq_C$ is not an isolated point in 
$\mathcal{LO}(F_2)$, but no method of constructing a sequence converging to $\preceq_C$ is given therein.  
Given an ordering $\preceq$ in $\mathcal{LO}(F_2)$, the known methods for constructing a sequence converging 
to $\preceq$ involve either approximation using the conjugates of $\preceq$, or approximation by modifying 
the ordering on the convex jumps in the Conradian soul of $\preceq$.  The results of this Appendix show that 
neither of these methods is sufficient for constructing a sequence of orderings converging to $\preceq_C$.
\end{rem}



\begin{small}

\vspace{0.1cm}


\noindent Adam Clay\\

\noindent Department of Mathematics, University of British Columbia, Vancouver, BC Canada V6T 1Z2\\ 

\noindent Email adress: aclay@math.ubc.ca\\

\noindent URL: http://www.math.ubc.ca/$\sim$aclay/ \\

\end{small}

\end{document}